\documentclass[12pt,oneside,reqno]{amsart}
\usepackage{lmodern}
\usepackage[T1]{fontenc}
\usepackage[english,activeacute]{babel}
\usepackage[applemac]{inputenc}
\usepackage{cite}
\usepackage{geometry}      		
\usepackage[normalem]{ulem}		
\usepackage{graphicx}			
\usepackage{amsmath}
\usepackage{amsfonts}
\usepackage{amssymb}
\usepackage{amsthm}
\usepackage{bbm}
\usepackage{MnSymbol}	
\usepackage{verbatim}
\usepackage{mathrsfs}
\usepackage[foot]{amsaddr}
\usepackage{color}
\usepackage{pgf}
\usepackage{tikz}

\usetikzlibrary{arrows,shapes,decorations,automata,backgrounds,petri,positioning}

\tikzset{main node/.style={circle,fill=blue!20,draw,minimum size=0.8cm,inner sep=0pt},}

\usepackage{layout}

\theoremstyle{definition}

\theoremstyle{remark}

\newcommand{\bs}[1]{\boldsymbol{#1}}

\def\>{\rangle}
\def\<{\langle}

\def\0{\bs{0}}
\def\1{\mathbbm{1}}

  \def\XXint#1#2#3{{\setbox0=\hbox{$#1{#2#3}{\int}$}
      \vcenter{\hbox{$#2#3$}}\kern-.47\wd0}}

\begin{document}


\title[A Feynman-Kac approach to a paper of Chung and Feller]
{A  Feynman-Kac approach to a paper of Chung and Feller on fluctuations in the coin-tossing game}
\author{\small F.~Alberto Gr\"unbaum$^1$}
\address{\scriptsize 
$^1$Department of Mathematics, University of California, Berkeley, CA 94720, USA}

\subjclass[2010]{60J10,60J65,33645,81Q30}

\keywords{Brownian motion, coin-tossing game, Feynman-Kac formula, Legendre polynomials} 

\date{}					

\begin{abstract}
A classical result of K.~L. Chung and W.~Feller deals with the partial sums
$S_k$ arising in a fair coin-tossing game. If $N_n$ is the number of ``positive'' terms among $S_1,S_2,\dots, S_n$ then the quantity 
$P(N_{2n} = 2r)$ takes an elegant form. We lift the restriction on an even number of tosses and give a simple expression for
$P(N_{2n+1} = r)$, $r = 0,1,2,\dots,2n+1$. We get to this result by adapting the Feynman--Kac methodology.
\end{abstract}

\maketitle

\tableofcontents

\section{Introduction}
\label{sec:INTRO}

One of the most surprising results about one dimensional Brownian motion goes back to
the remarkable work of P.~L\'evy, see \cite{L}, who discovered the famous
arcsine law for the proportion of time that the particle spends to the right of the origin. The arcsine law reappeared shortly later in a general result of P.~Erd\H{o}s and M.~Kac for partial sums of independent random variables, see \cite{EK}. The arcsine law has been revisited in many publication, see for instance \cite{P1} and references in this paper.

This amazing result of P.~L\'evy was reproved around 1950 by M.~Kac, see \cite{K}, using
what was at the time the recently developed Feynman--Kac method of connecting integration on
function space with solutions of diffusion equations, beyond the well-known
case of the heat equation with a zero potential.
Feynman had proposed expressing transition amplitudes in quantum mechanics in terms of a Lagrangian, see \cite{Fe1,Fe2} and M. Kac who by then had done extensive work on functionals of Brownian motion adapted these heuristic ideas and made them rigorous in the case of classical probability. For an instructive view of the relation between the work of R. Feynman and that of M. Kac the reader may refer to \cite{K1,K2,S}.
Seminal efforts to use a Lagrangian in quantum mechanics go back to P. Dirac, see \cite{D}.

\bigskip

Of about the same vintage is the first detailed study 
of discrete analogs of the arcsine law for Brownian motion
undertaken in the paper by K.~Chung
and W.~Feller, see \cite{CF}, dealing with fluctuations in the coin-tossing game. In this case ``the constant reference to ties can render the exposition clumsy'' (quoting from \cite{F}) and a good way out of this difficulty was given in \cite{CF} as we recall now.

The discrete case requires a careful agreement of what one means by
the partial sum $S_k$ being positive at time $k$ if one is to get ``elegant 
results'' (quoting from \cite{CF}). One says that $S_k$ is positive if 
$S_k > 0$ or in case that $S_k=0$ one has $S_{k-1} > 0$. If $N_n$ denotes the number 
of ``positive'' terms among $S_1,S_2,\dots,S_n$ one finds that
\[
P(N_{2n} = 2r)  = \frac{\binom{2r}{r}  \binom{2n-2r}{n-r}}{2^{2n}}\,.
\]
Using the notation $u_{2k}=\frac{\binom{2k}{k}}{2^{2k}}$ we have 
\[
P(N_{2n} = 2r)  = u_{2r} u_{2n-2r}\,.
\]
The purpose of this paper is to give a simple expression for the quantities

\begin{equation}\label{odd1}
P(N_{2n+1}=r)
\end{equation}
which we have failed to find in the literature. One should remark that
in the first edition of the celebrated book by W.~Feller, see \cite{F}, footnote~4, page~252,  one finds the comment that the method in \cite{CF} leads to an expression for these quantities. This remark does not appear in later editions of this book. The recent book, \cite{McK}, page~78, poses this as an exercise.
W.~Feller used to refer to the act of publishing results that are already known as ``educating yourself in public''. I hope that this is not one such case, or at least that the expressions given here may be simpler than other known ones. Our result reads as follows
\begin{equation}\label{modd}
P(N_{2n+1}=2r)=u_{2r} u_{2n+2-2r} \frac{n-r+1}{n+1}, r=0,1,2,\dots,n
\end{equation}

and 

\begin{equation}\label{nodd}
P(N_{2n+1}=2r-1)=u_{2r} u_{2n+2-2r} \frac{r}{n+1}, r=1,2,\dots,n+1\,.
\end{equation}

\bigskip

Given the success of the Feynman--Kac method in reproving the result of P.~L\'evy for Brownian motion it is a bit surprising that this method has not been exploited to get the Chung--Feller result for coin tossing.

\bigskip

There is a nice paper, see \cite{C},
where the author discusses a discrete version of the Feynman--Kac formula.
This paper considers the general case when a function $p(n,x)$ is defined by the expression
\begin{equation}\label{expec}
	p(n,x) = p(n,x,q) = E_x\left( q^{\sum_{k=0}^n f(S_k)} g(S_n)\right)\,.
\end{equation}

\noindent
with $f(x),g(x)$ arbitrary functions. The subindex $x$ indicates that the walk $S_k$ starts at position $x$, and $q$ is an undeterminate. The analogy with M. Kac's formulation is clear from \eqref{expectedv} below.

Very much in the spirit of \cite{K}, ones sees in \cite{C} that this function $p(n,x,q)$ satisfies a partial difference equation, in the variables $(n,x)$, and then the function 
\begin{equation}\label{double}
\psi(x,q,z) = \sum_{n=0}^{\infty} p(n,x,q)z^n,
\end{equation}

\noindent
which is a double generating function of the distribution of interest,
satisfies an appropriate ordinary difference equation, which one can solve.

This allows E. Csaki to discuss in -principle- see
Example~$3.1.3$ in \cite{C}, the distribution of the number of terms where $S_k$, $k=0,1,\dots,n$ is non-negative. The explicit expression for $\psi(0,q,z)$ 
going with this example will be given here. To get the distribution in question one still needs to invert a double generating function. The expression for $\psi(0,q,z)$ is essentially given in \cite{C}, but the inversion is not pursued in that paper or here.

\bigskip

One can see that, as Chung and Feller already knew,
this random variable (the number of terms when $S_k$ is non-negative) has a different distribution than the one considered in \cite{CF}.
In the case of Brownian motion these differences go
away.

To summarize, Example~$3.1.3$ in \cite{C} deals with 
the random variable
\[
q^{\sum_{k \le n} {1\negthinspace\negthinspace 1}\{S_k \ge 0\}}
\]
which fits with the general method in that paper for an appropriate choice of the functions $f(x),g(x)$, namely $f$ is the characteristic function of the non-negative integers and $g$ is identically one. 

Given the way in which ``positive $S_i$'' was defined in \cite{CF} one cannot write the random variable $q^{N_n}$ in the form $(1.6)$ of \cite{C}.

\bigskip

It is however still possible to consider 
\[
	p(n,x)= p(n,x,q) = E_{x} (q^{N_n})
\]
and derive a partial difference equation for $p(n,x)$.

Once again, one defines a function $\psi(x,q,z)$, as in
\eqref{double}.
This function satisfies, once again, an ordinary difference equation which can be solved using the same method given in \cite{C}. In other words, we
see that the method in \cite{C} can be applied to some cases beyond those written as in expression $(1.6)$, page~64, \cite{C}.

\section{Comparing with the original Feynman-Kac method}

The discussion above is a perfect analog of the treatment in \cite{K}
where one starts with a partial differential equation in the variables $(x,t)$ satisfied
by a function $p(t,x)$ defined by an expectation as in \eqref{expec}.

More explicitly, Kac
considers the function 
\begin{equation}\label{expectedv}
	u(t,x) = E_{x} (e^{\int_0^t V(b(s,w))ds}f(b(t,w)))
\end{equation}

\noindent
Here $b(s,w)$ stands for one dimensional Brownian motiom.
Now the Feynman-Kac approach says that $u(x,t)$, defined above, solves the initial value problem for a partial differential equation
\begin{align*}
u_t &= \frac {1}{2} u''(t,x) + V(x)u, \text{ with} \\
u(0,x) &= f. 
\end{align*}

\bigskip

Kac considers the partial differential equation above for the family of functions 
$u(t,x,\beta)$ that come about when one takes the expected value of the functionals of Brownian motion given by

\[
e^{-\beta \int_0^t V(b(s,w)ds}
\]

\noindent
with $V(x)$ the characteristic function of the set ${x \ge \0}$ and where $f$ has been chosen to be identically one in 
\eqref{expectedv}. 
This introduces a Laplace transform (in the variable $\beta)$ of the distribution function for the random variable of interest.

The classical way of solving an initial value problem for a partial differential equation as the one above is to take its Laplace transform in the time variable $t$

\[
	{\hat u} (x,\alpha,\beta) =  \int_0^{\infty} e^{-\alpha t} u(t,x,\beta) dt
\]

This function of $x$ solves now an inhomogeneous linear second order ordinary differential equation. Its solution is given by finding, and then using, the appropriate Green function, i.e. by finding two linearly independent solutions of the homogeneous equation and matching their product properly at the discontinuity. M. Kac succeeded in carrying up all these steps and by evaluating ${\hat u}(x,\alpha,\beta)$ at $x=0$ he obtained the expression

\[
{\hat u}(0,\alpha,\beta) = \frac {1}{\sqrt{\alpha(\alpha+\beta)}}
\]
which is nicely given already as a double Laplace transform 
\[
\frac {1}{\sqrt{\alpha(\alpha+\beta)}} = \int_0^{\infty} e^{-\alpha t} \left( \frac {1}{\pi} \int_0^t \frac {e^{-\beta s}}{\sqrt{(t-s)s}} ds\right)dt
\]
This makes the inversion automatic and yields the arc-sine law. See \cite{K}, page~192, and for more details, see \cite{IM}, pages 54--57, 
or \cite{S1}, pages 230--233.

Now that we have reviewed the original Feynman-Kac approach we see that
in the discrete case E. Csaki replaces $e^{-\beta}$ by $q$ and $e^{-\alpha}$ by $z$ and starting from Kac's formula, \eqref{expectedv}, he replaces the functions $(V,f)$ by $(f,g)$ to arrive at the expression (1.6), page 64 in \cite{C}, i.e. \eqref{expec} above.

\bigskip

The only difference between \cite{K} and what we will do in this paper is that
the inversion of the double generating function in the discrete case will not be so simple as the case recalled above.

\bigskip

The main result of this paper is the expression for the function $\psi(0,q,z)$, see \eqref{double}.
The function $\psi(0,q,z)$, as a function of $z$, can be written as the sum of its even part

\[
\psi_{even}(0,q,z) = \sum_{n=0}^{\infty} p(2n,0,q) z^{2n}
\]

and its odd part

\begin{equation}\label{oddexpress}
\psi_{odd}(0,q,z) = \sum_{n=0}^{\infty} p(2n+1,0,q) z^{2n+1}
\end{equation}

\bigskip

The even part was already given in \cite{CF,F} and the odd part is new, see \eqref{express}.
Getting our hands on its coefficients $p(2n+1,0,q)$ will yield \eqref{modd} and \eqref{nodd}.

\noindent
In both cases, even or odd $m$, the generating function for $N_{m}$ is given by 

\[
	p(m,0,q)=\sum_{q=0}^{m} P(N_{m}=j) q^j
\]


\section{The problem of Chung and Feller}

We first supplement the results of \cite{C} by giving the value of $\psi(0,q,z)$
corresponding to Example~$3.1.3$.

One gets as the solution to the appropriate ordinary difference equation, evaluated at $x=0$ the expression
\[
\psi(0,q,z)=\frac {q}{1-qz} + \frac {\sqrt{1-q^2-z^2}(\sqrt{1-z^2} + 1-z) - (1-qz)\sqrt{1-z^2} + (1-q(1-z))z - 1}{2(z-1)z(qz-1)}\,.
\]

This is the double generating function whose inversion would give the distribution of the number of $S_k \geq 0$, $k=0,1,\dots,n$.

\bigskip

\bigskip

Now we move to the situation relevant to the Chung--Feller result, \cite{CF}.

\bigskip
The partial difference equation satisfied by $p(n,x)=p(n,x,q)$ in this case is given as follows

\bigskip
if $x > 0$ \qquad $p(n,x) = q/2(p(n-1,x+1)+p(n-1,x-1))$

\medskip
if $x = 0$ \qquad $p(n,x) = q/2\,p(n-1,x+1) + 1/2\,p(n-1,x-1)$

\medskip
if $x < 0$ \qquad $p(n,x) = \frac {1}{2} (p(n-1,x+1) + p(n-1,x+1))$

\bigskip

One can use the general analytical method put forward in \cite{C}, to solve for

\[
\psi(x,q,z) = \sum_{n=0}^{\infty} p(n,x,q)z^n,
\]

The approach in \cite{C} gives three different expressions for $\psi(x,q,z)$
depending on $x$ being positive, zero or negative, namely

\bigskip
\[
\begin{aligned}
x > 0 &\quad\psi_{\mbox{posit}}(x,q,z) = e(z) \left( \frac {1-\sqrt{1-z^2q^2}}{zq}\right)^x + 
\frac {1}{1-zq} \\
x = 0 &\quad\psi_{\mbox{zero}}(x,q,z) = c(z) \left( \frac {1-\sqrt{1-z^2q}}{zq}\right)^x + 
\frac {2}{2-zq-z} \\
x < 0 &\quad\psi_{\mbox{negat}}(x,q,z) = d(z) \left( \frac {1+\sqrt{1-z^2}}{z}\right)^x + 
\frac {1}{1+z}
\end{aligned}
\]

with $z$ dependent constants $e(z),c(z),d(z)$ to be determined by imposing the 
conditions

\bigskip
\[
\begin{aligned}
\psi_{\mbox{posit}}(1) &= 1 + \frac {qz}{2} (\psi_{\mbox{posit}}(2) + \psi_{\mbox{zero}}(0)) \\
\psi_{\mbox{zero}}(0) &= 1 + \frac {qz}{2} \psi_{\mbox{posit}}(1) + \frac {z}{2} \psi_{\mbox{negat}}(-1) \\
\psi_{\mbox{negat}}(-1) &= 1 + \frac {z}{2} (\psi_{\mbox{zero}}(0) + \psi_{\mbox{negat}}(-2))
\end{aligned}
\]

\bigskip

Eventually we are only interested in the value of $\psi(x,q,z)$ when $x=0$ since that is the initial value for our walk.
Things simplify considerably when one sets $x=0$ and one gets a (relatively) simple analytical expression for $\psi(0,q,z)$ displayed below in \eqref{nimer} and \eqref{denim}.

\bigskip

After some minor simplifications this can be written as the ratio of the expressions

\begin{equation}\label{nimer}
-z(q+1)\sqrt{1-q^2z^2}((qz^2-1)(2z^2+\sqrt{1-z^2}(z^2-1)-1) \\
\quad + (1-z^2)\sqrt{1-q^2z^2}(2\sqrt{1-z^2}-z^2+2))
\end{equation}

\bigskip

and

\begin{equation}\label{denim}
(1-z^2)(1-q^2z^2)(((q^2+1)z^2-2)(2\sqrt{1-z^2}-z^2+2) \\
\quad + \sqrt{1-q^2z^2}(4z^2+2\sqrt{1-z^2}(z^2-2)-4))
\end{equation}

Starting from \eqref{nimer} and \eqref{denim} the even part of $\psi(0,q,z)$ i.e. $$\frac{\psi(0,q,z)+\psi(0,q,-z)}{2}$$
is given by
\[
\frac{1}{\sqrt{1-z^2}\sqrt{1-q^2 z^2}}
\]

\bigskip

\bigskip
 
The result about the even part is already contained in \cite{CF,F}. By using
the generating function of the Legendre polynomials $P_n$ (with $P_n(1)=1$) we get a partial inversion
of the double generating function for an even number of coin-tosses, namely

\[
	p(2n,0) = p(2n,0,q)= \sum_{k=0}^{n} u_{2k} u_{2n-2k} q^{2k}=P_n( \frac{q+q^-1}{2}) q^n\,.
\]

This connection of the Legendre polynomials and the coin-tossing game is  
known, see \cite{R}, formula(10), page 248.
Since the actual probabilities $u_{2k} u_{2n-2k}$ are known, this is a case where the double generating function can be fully inverted.
It will be convenient to introduce, for later use, the notation

\begin{equation}\label{theas}
A(n,q)=P_n\left( \frac{q+q^-1}{2}\right) q^n. 
\end{equation}

\bigskip

So far we have used the Feynman-Kac method given in \cite{C} to reprove the results in \cite{CF}.


\bigskip

\bigskip

\noindent
From now on we concentrate the odd part of $\psi(0,q,z)$

$$\psi_{odd}(0,q,z) = \frac{\psi(0,q,z)-\psi(0,q,-z)}{2}$$

\noindent
which is given by the ratio of the expressions
\begin{equation}\label{num1}
\sqrt{1-z^2}\sqrt{1-q^2z^2}(qz^2+1)-z^2(q^2(z^2-1)-1)-1 
\end{equation}
and
\begin{equation}\label{denom1}
(1-z^2)(1-q^2z^2)(q+1)z,
\end{equation}

\noindent
This ratio can be conveniently rewritten in the form

\[
\psi_{odd}(0,q,z)=\frac {qz^2+1-\sqrt{1-z^2} \sqrt{1-z^2q^2}}{(q+1)z\sqrt{1-z^2} \sqrt{1-z^2q^2}}
\]

or in the alternative form

\begin{equation}\label{express}
\psi_{odd}(0,q,z)=\frac {\frac {1}{\sqrt{1-z^2} \sqrt{1-z^2q^2}} - 1}{(q+1)z} + \frac {qz}{\sqrt{1-z^2} \sqrt{1-z^2q^2} (q+1)}
\end{equation}

\noindent
In the expression above
we look at the part involving the inverse of the product of two square roots.
If one puts

\[
\frac {1}{\sqrt{1-z^2} \sqrt{1-z^2q^2}} = \sum_{i=0}^{\infty} v_iz^{2i}
\]

\noindent
and uses $v_0=1$,
it is easy to see that $\psi_{odd}(0,q,z)$, i.e. \eqref{express}, is given as follows

\[
\sum_{i=0}^{\infty} \frac {v_iq + v_{i+1}}{q+1} z^{2i+1}
\]

This implies, see \eqref{oddexpress}, that

\[
\begin{aligned}
p(2n+1,0,q) &= \frac {A_n(q)q + A_{n+1}(q)}{q+1} 
\end{aligned}
\]

It is easy to see that if one puts

\begin{equation}
	A_n(q)=\sum_{i=0}^{n} w[2i,2n] q^{2i}
\end{equation}

and thus

\begin{equation}
	A_{n+1}(q)=\sum_{i=0}^{n+1} w[2i,2n+2] q^{2i}
\end{equation}

then the ratio

\begin{equation}
p(2n+1,0,q) = \frac {A_n(q)q + A_{n+1}(q)}{q+1} 
\end{equation}

has the expansion

\begin{equation}
\sum_{i=0}^{2n+1} x_i q^i
\end{equation}

with the coefficients $x_i$ given as follows

\begin{equation}
	x_{2i} = \sum_{j=0}^{i} w[2j,2n+2] - \sum_{j=0}^{i-1} w[2j,2n]
\end{equation}

and

\begin{equation}
	x_{2i+1}=\sum_{j=0}^{i} w[2j,2n] - \sum_{j=0}^{i} w[2j,2n+2]
\end{equation}

Finally using the fact that $$w[2j,2n]=u_{2j} u_{2n-2j}$$
and recalling that

$$x_i = P(N_{2n+1}=i)$$
one can prove \eqref{modd}   and \eqref{nodd}, for instance by induction.

\section{Surmising the new results}
Now I move to an account of how these results were conjectured, before the proof given above was found.

With the expression for
$\psi_{odd}(0,q,z)$ at hand, see \eqref{num1} and \eqref{denom1}, I could compute, with the assistance of the symbol manipulator Maxima, as many of the coefficients $p(2n+1,0,q)$ as necessary.
These are polynomials in $q$ which allow one to give, in principle, the discrete
analog of the arcsine law for any number of tosses $2n+1$.

At this point I did a bit of experimentation. Plotting the values of 
$P(N_{2n+1}=j, j=0,1,2,....,2n+1)$ was not particularly useful. However plotting the cummulative distribution functions for three consecutive values of time namely $2n,2n+1,2n+2$ revealed an interesting pattern which in turn lead me to conjecture the expressions given in
\eqref{modd} and \eqref{nodd}. Here it was very important that the results for an even number of tosses were already known.

Having checked my claims
\eqref{modd} and \eqref{nodd}
in a number of ways, I showed them 
to my colleague Jim Pitman. He told me he had never seen these explicit expressions before, but in a few days he found a probabilistic proof of these claims, and then proceeded to connect these results with some of his and Marc Yor's work from several years ago. He also showed me how these claims are the exact discrete versions of some results of P. L\'evy. I am extremely thankful to him for explaining to me how these new results connect with very deep and beautiful work. Some relevant references are \cite{JPY,P2,PY}. I hope that sometime soon he will write some of this up. A minimal indication is given below.

\bigskip

In \cite{PY} the authors recall as (4.m)  
 formula (51) of 
  P. L\'evy \cite{L} 
that if

$$A_+(t):= \int_0^1 1(B_s > 0 ) ds$$  then
$$
P( t ^{-1} A_+(t) \in du , B_t > 0 ) = u P( A_+(1) \in du ) = \frac{du} {\pi} u^{1/2} (1-u)^{-1/2} \qquad ( 0 < u < 1 ).
$$.

In combination with the arcsine law this gives (4.n) in \cite{PY}, see below.
J. Pitman pointed out to me \cite{P4} that 
\eqref{modd} and \eqref{nodd}, properly rewritten 
amount to a combination of the Chung-Feller
discrete arcsine formula with the simple conditional probability 
\begin{equation}
\label{cond1}
P( S_{2n - 1} > 0 \,|\, N_{2n} = 2 r)  = \frac{ r } {n }  \qquad ( 0 \le r \le n ).
\end{equation}
This is an exact discrete analog of the consequence of L\'evy's Brownian formula,  observed in
(4.n) of \cite{PY}, that
\begin{equation}
\label{condbm}
P( B_t > 0  \,|\, A_+(t) = a )  = \frac{ a } {t } \qquad ( 0 \le a \le t )  .
\end{equation}

\bigskip

\bigskip
I close this section by giving  a few more expressions for $p(2n+1,0,q)$ in terms of the quantities introduced in \eqref{theas}. One such formula is

\[
\begin{aligned}
p(2n+1,0,q)=A(n+1,q) + \frac{1-q}{2(n+1)} \frac{d}{dq} A(n+1,q)
\end{aligned}
\]

Using now the differentiation formula for the Legendre polynomials, this can be
expressed as

\[
\begin{aligned}
&p(2n+1,0) = \frac {q}{q+1} \frac {n+2}{2n+3} A(n,q) \\
&\quad + \frac {q+1}{2q} A(n+1,q) - \frac {1}{q(q+1)} \frac {n+2}{2n+3} A(n+2,q)
\end{aligned}
\]

which can be rewritten in several ways, such as

\[
\begin{aligned}
	&p(2n+1,0) = \frac {((2n+3)(q^2+1)+(2n+2)q)q A_n(q)+2(n+2) A_{n+2}(q)}{(1+q+q^2+q^3)(2n+3)} 
\end{aligned}
\]

Finally, one more expression is

\bigskip
\[
\begin{aligned}
p(2n+1,0,q)= \frac {A_{n+1}(q) - q^2A_n(q)}{1-q^2} + q \frac {(A_n(q) - A_{n+1}(q))}{1-q^2}
\end{aligned}
\]

where we have decomposed $p(2n+1,0,q)$ into its even and odd parts as a function of $q$.



\section{Concluding remarks}

My motivation to revisit the classical results on the coin-tossing game was
an (still in progress) effort to consider this topic in the case of quantum walks. This is joint work  with L. Velazquez and others. 

Some readers may be intrigued by the appearance of the Legendre polynomials in a probability problem.
Richard Askey pointed 
out to me that the Legendre polynomials had actually surfaced for the first time in connection with some problem in probability theory.
He recalls that Arthur Erdelyi told him once that this occurred in some work by J. L. Lagrange. Indeed in \cite{L} , Lagrange considers the problem of repeated observations of an unknown quantity. A certain number of times, denoted by $a$, one measures the exact value, while for $b$ times the error is $+1$ and for $b$ times the error is $-1$. He is interested in the probability of getting the exact result by using "le milieu" -the mean- of these repated measurements.

In the process of solving this, he needs to find the power series expansion of 
the expression

\[
   \frac {1}{\sqrt{1-2az+(a^2-4b^2)z^2}}
\]

\noindent
in powers of $z$.

\bigskip

The three term recursion for the coefficients in this expansion is explicitly written down and considered well known
by Lagrange. The case $a^2-4b^2=1$ gives the Legendre polynomials in the variable $a$. The work of Lagrange took place in the period 1770-1773, and predates the work of Legendre and Laplace.

\bigskip

A final comment about the Feynman-Kac approach, whose validity goes beyond one
dimensional Brownian motion. One should not assume that solving the relevant ordinary differential equation is always an easy matter. For an example where I have not been able to solve
the appropriate inhomogeneous second order ordinary differential equation see
 \cite{GM}. The case in \cite{K} and the one here where the equation can be solved and then the double Laplace transform or the double generating function can be explicitly inverted are rather exceptional.

\end{document}